\documentclass[12pt, a4paper]{article}
\usepackage[left=1in,right=1in,top=1.0in,bottom=1in]{geometry}

\usepackage{amsmath,amssymb,amsthm,amsfonts,afterpage,hyperref}
\usepackage{amscd,subeqnarray}

\usepackage{multirow}

\usepackage{graphicx}
\usepackage{subcaption}
\usepackage{authblk}
\usepackage{tikz}

\usepackage{array}
\usepackage{wrapfig}
\usepackage{color}
\usepackage{ragged2e}
\usepackage{stmaryrd}
\usepackage{siunitx}
\usepackage{multirow}
\usepackage{xcolor,cancel}
\usepackage{boldfonts}
\usepackage{symbols}
\usepackage{footmisc}
\usepackage{float}
\usepackage{setspace}
\usepackage{bm}
\usepackage{booktabs}
\usepackage{tikz, xcolor}
\usepackage{soul}

\usepackage{tabularx} 
\usepackage{rotating}

\usepackage{algorithm,algorithmic} 
\usepackage[linesnumbered,ruled,vlined,algo2e]{algorithm2e}
\SetKwInput{KwInput}{Input}                
\SetKwInput{KwOutput}{Output}

\usepackage[sort, numbers]{natbib}
\setcitestyle{square}

\usetikzlibrary{shapes,arrows, positioning,calc}	
\usetikzlibrary{arrows,snakes,backgrounds}

\usepackage{enumitem} 
\usepackage{ragged2e} 

\newcommand{\inputlayer}{\textbf{input layer}}
\newcommand{\hiddenlayer}{\textbf{hidden layer}}
\newcommand{\outputlayer}{\textbf{output layer}}
\newcommand{\weightvector}{\textbf{weight vector}}
\newcommand{\bias}{\textbf{bias}}
\newcommand{\activationfunction}{\textbf{activation function}}
\newcommand{\backpropagation}{\textbf{backpropagation algorithm}}
\newcommand{\lossfunction}{\textbf{loss function}}
\newcommand{\differentialcost}{\textbf{differential cost}}
\newcommand{\boundarycost}{\textbf{boundary cost}}
\newcommand{\gradientdescent}{\textbf{gradient descent methods}}
\newcommand{\semismoothNewton}{\textbf{semi-smooth Newton's methods}}

\DeclareMathAlphabet{\mathpzc}{OT1}{pzc}{m}{it}

\theoremstyle{definition}
\newtheorem{problem}{Problem}
\newtheorem{weakform}{Weak Formulation}

\theoremstyle{plain}
\newtheorem{theorem}{Theorem}

\theoremstyle{remark}
\newtheorem{remark}{Remark}

\providecommand{\keywords}[1]
{
  \small	
  \textbf{\textit{Keywords---}} #1
}

\usepackage{array}
\newcolumntype{C}{>{\raggedright\arraybackslash}p{2.5cm}}

\title{\textsf{Numerical Approximation of Electrohydrodynamics Model: A Comparative Study of PINNs and FEM}}

\author[1]{Mara Martinez \thanks{mmcguire4@islander.tamucc.edu}} 
 \author[1]{B. Veena S. N. Rao \thanks{bv.rao@tamucc.edu}} 
\author[1]{S. M. Mallikarjunaiah\thanks{m.muddamallappa@tamucc.edu}\thanks{corresponding author}} 
\affil[1]{Department of Mathematics \& Statistics,
Texas A\&M University-Corpus Christi, 
6300 Ocean Dr., 
Corpus Christi, Texas 78412-5825, USA}
\date{}

\begin{document}

\maketitle

\begin{abstract}
The accurate representation of numerous physical, chemical, and biological processes relies heavily on differential equations (DEs), particularly nonlinear differential equations (NDEs). While understanding these complex systems necessitates obtaining solutions to their governing equations, the derivation of precise approximations for NDEs remains a formidable task in computational mathematics. Although established techniques such as the finite element method (FEM) have long been foundational, remarkable promise for approximating continuous functions with high efficacy has recently been demonstrated by advancements in physics-informed deep-learning feedforward neural networks. In this work, a novel application of PINNs is presented for the approximation of the challenging Electrohydrodynamic (EHD) problem. A specific $L^2$-type \textit{total loss function} is employed, notably without reliance on any prior knowledge of the exact solution. A comprehensive comparative study is conducted, juxtaposing the approximation capabilities of the proposed neural network with those of the conventional FEM. The PINN training regimen is composed of two critical steps: forward propagation for adjustments to gradient and curvature, and backpropagation for the refinement of hyperparameters. The critical challenge of identifying optimal neural network architectures and hyperparameter configurations for efficient optimization is meticulously investigated. Excellent performance is shown to be delivered by the neural network even with a limited training dataset. Simultaneously, it is demonstrated that the accuracy of the FEM can be substantially enhanced through the judicious selection of smaller mesh sizes.
\end{abstract}

\noindent \keywords{Unsupervised machine learning, Activation function, Optimizer,  Darcy-Brinkman-Forchheimer model, Network sensitivity, Finite element method }

\section{Introduction}
The mathematical modeling of complex phenomena across engineering and the natural sciences frequently culminates in the formulation of nonlinear differential equations. These equations provide a framework for faithfully representing the intricate behaviors of systems ranging from fluid dynamics \cite{vasilyeva2023,SMMDDB2023} and solid mechanics \cite{lee2022finite,HyunMSM_MMS2022,yoon2022preferential,yoon2021finite,Mallikarjunaiah2015,yoon2021quasi,ferguson2015numerical,rao2022mathematical,MVSMM2023,gou2015modeling,Gou2023a,Gou2023b,manohar2024hp} to abstract concepts like nonlocality in boundary value problems \cite{muddamallappa2015two,walton2016plane}. At their core, these models are constructed by deriving governing equations from fundamental physical laws that describe the spatial and temporal variation of key quantities. This process gives rise to a central question in applied mathematics: given a well-posed model, which guarantees the existence and uniqueness of a solution, how can this unknown solution function be effectively constructed? For a meaningful interpretation of a system's behavior, it is imperative to obtain either a closed-form analytical solution or a high-fidelity numerical approximation.

However, the analytical treatment of NDEs arising from real-world applications is often intractable due to inherent nonlinearities, geometric complexities of the domain, or sophisticated boundary conditions. This necessitates a turn toward numerical and approximate methods. The fundamental strategy of these techniques is to transform the continuous problem defined by the NDE into a discrete one, thereby yielding a system of algebraic equations that can be solved computationally. The solution to this system provides an approximation of the true analytical solution at a finite set of points across the domain. Within the suite of powerful numerical techniques, the FEM has achieved widespread adoption. In FEM, the numerical solution is expressed as a linear combination of piecewise polynomial basis functions, with unknown scalar coefficients representing the solution's nodal values. Other established and extensively employed methods for solving DEs include spectral methods, finite difference methods, and finite volume methods.

In a recent paradigm shift, the confluence of machine learning and scientific computing has introduced artificial neural networks (ANNs) as a compelling, mesh-free alternative for function approximation \cite{lee1990neural,meade1994solution,meade1994numerical,lagaris2000neural}. Specifically, Physics-Informed Neural Networks (PINNs) represent a potent methodology where the governing equations are directly integrated into the network's learning algorithm. This is achieved by designing a composite loss function that penalizes the network's predictions for deviations from not only available data points but also the underlying physical constraints imposed by the differential equation and its boundary conditions. This physics-informed regularization guides the optimization process, enabling PINNs to converge to accurate solutions even with sparse or noisy data. A key operational advantage is their mesh-free nature, which circumvents the often laborious and computationally intensive process of grid generation required by methods like FEM. Furthermore, all necessary derivatives are computed with machine precision via automatic differentiation, precluding the accumulation of discretization errors. While challenges related to the optimization landscape and the balancing of loss terms remain active areas of research, the ability of PINNs to synergistically leverage both data and first principles makes them an exceptionally promising tool for tackling complex DEs.

Inspired by the architecture of biological nervous systems, ANNs are machine learning algorithms whose capacity for function approximation is theoretically guaranteed by the Universal Approximation Theorem (UAT) \cite{cybenko1989approximation,hornik1989multilayer}. The UAT posits that a feedforward network with a single hidden layer containing a sufficient number of neurons can approximate any continuous function to an arbitrary degree of accuracy. This theoretical foundation is complemented by the practical advantage of automatic differentiation, a technique that allows for the exact and computationally efficient calculation of derivatives of any order for the network's output with respect to its inputs \cite{baydin2018automatic,rall1981automatic}. Despite these powerful features, many early applications of ANNs to differential equations were limited. For instance, some approaches imposed boundary conditions externally by structuring the network's output as part of a trial solution, while others relied on the exact solution as a training target, limiting their utility \cite{dockhorn2019discussion}. Moreover, a common critique of prior work is the treatment of ANNs as ``black-box'' optimizers, often without a rigorous justification for the chosen loss function or a systematic analysis of the hyperparameter sensitivities that govern solution accuracy \cite{SMM2023,venkatachalapathy2022,venkatachalapathy2023,piscopo2019solving,michoski2020solving,Mara2023b,hussain2025machine,zafar2025optimizing,haider2025machine}.

The present study extends previous research \cite{SMM2023,venkatachalapathy2022,venkatachalapathy2023,Mara2023b} by introducing a deep learning architecture expressly developed for approximating complex NDEs, and it addresses the aforementioned gaps in the literature. A primary contribution of this work is a systematic investigation into the optimal hyperparameter configurations required for robust and accurate function approximation. We introduce a mean-squared error metric to formalize the performance evaluation of the network. The total loss function is formulated as a composite of a differential cost and a boundary loss, critically enabling the training process to proceed with only an approximate solution, a significant advantage for practical problems where analytical solutions are unavailable. A sparse training dataset is constructed from the computational domain, and the network autonomously learns the solution's behavior, adjusting for gradient and curvature through a fixed number of training iterations. To contextualize the performance of our deep learning model, a rigorous comparative analysis against a continuous Galerkin FEM is performed, evaluating both solution accuracy and computational cost. Our findings demonstrate that the network-based approach can achieve high accuracy without the need for re-meshing, a common requirement for FEM, and is readily extensible to higher-dimensional problems.

The remainder of this manuscript is organized as follows. Section~\ref{background} establishes the theoretical background by reviewing relevant function spaces and formally defining the one-dimensional Electrohydrodynamic problem. Section~\ref{sec:ANN_architecture_training} details the proposed neural network architecture and the associated training algorithm. The numerical experiments applying this network to the nonlinear boundary value problem are presented in Section~\ref{num_exp}. For a direct comparison, the finite element discretization of the problem is outlined in Section~\ref{sec:numerical_discretization}. A comprehensive comparative analysis of the two methods is conducted in Section~\ref{comparision}. Finally, Section~\ref{conclusion} provides a summary of the key conclusions and suggests directions for future research.

\section{Mathematical Model}\label{background}
\begin{figure}[H]
\centering
\includegraphics[width=0.4\textwidth]{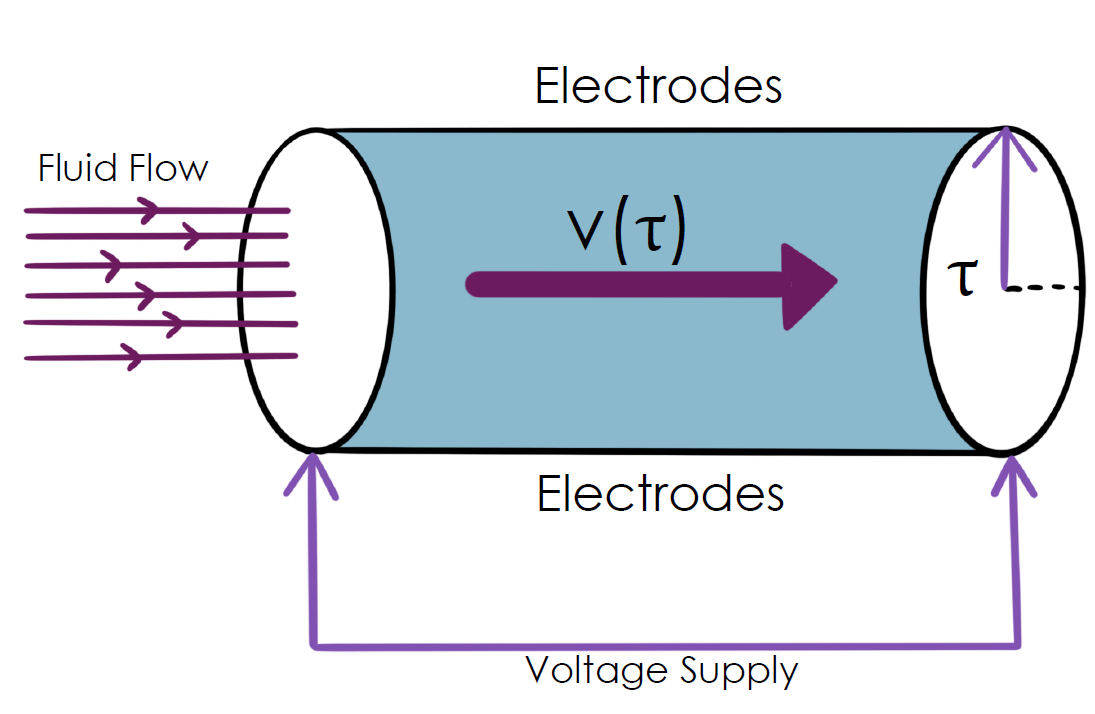}
\caption{Configuration showing EHD fluid flow in a circular cylindrical conduit.  }
\label{fig_disp_1a}
\end{figure}

This investigation presents a numerical framework for solving the complex mathematical model that governs the Electrohydrodynamic (EHD) flow of a fluid within a circular cylindrical conduit. The central focus is the high-fidelity approximation of this system, which is of significant practical importance for understanding phenomena such as ``ion draught'' configurations, depicted schematically in Figure~\ref{fig_disp_1a}. The EHD flow of a fluid in a circular cylindrical conduit is a system of significant practical importance for modeling phenomena like ``ion draught,'' which presents a challenging nonlinear problem.  A defining feature of this study is the rigorous benchmarking of our proposed methodology, which is based on the architecture of a PINN, against the well-established FEM. This direct comparison is designed to elucidate the relative performance, accuracy, and computational efficiency of the PINN approach for this class of intricate multiphysics problems.

The physical system is mathematically described by the following second-order nonlinear ordinary differential equation:
\begin{equation}\label{EHD_DE}
\dfrac{d^2 v}{d \tau^2} + \dfrac{1}{\tau} \, \dfrac{d v}{d \tau} + Ha^2 \left( 1 - \dfrac{v}{(1 - \alpha \, v)} \right) = 0, \quad 0 < \tau < 1,
\end{equation}
which is subject to the boundary conditions:
\begin{equation}\label{EHD_Bcs}
v^{\prime}(0) = 0, \quad v(1) = 0.
\end{equation}
Within this formulation, the independent variable $\tau$ denotes the radial distance from the cylinder's central axis, while the dependent variable $v(\tau)$ represents the fluid velocity. The model is characterized by two key dimensionless parameters: the Hartman number, $Ha$, which quantifies the influence of the electric field, and the parameter $\alpha$, which governs the degree of the system's nonlinearity. The value of $\alpha$ is a function of several physical properties, including the fluid's pressure gradient, ion mobility, and the inlet current density. For a comprehensive derivation of this boundary value problem (BVP), the reader is referred to the foundational works of \cite{mckee1997calculation,beg2013chebyshev,hosseini2017numerical}.

The solution to this problem has been pursued through a diverse spectrum of numerical and semi-analytical techniques. Notable contributions to the literature include, but are not limited to, the methods presented in \cite{abukhaled2021fast,sabir2018electrohydrodynamic,tiwari2018solution,hosseini2017numerical,shiralashetti2016haar,moghtadaei2015spectral,ghasemi2014electrohydrodynamic,beg2013chebyshev,hosseini2017numerical,hasankhani2016application,pandey2012semi,mastroberardino2011homotopy,paullet1999solutions,mckee1997calculation,sahlan2022lucas,thirumalai2019spectral,beg2013dtm,mastroberardino2011homotopy,khan2012approximate,roul2019new}.

\begin{remark}
It is important to note that for a non-trivial solution to exist, the Hartman number must be non-zero ($Ha^2 \neq 0$). In the case where $Ha^2 = 0$, the BVP \eqref{EHD_DE}-\eqref{EHD_Bcs} admits only the trivial solution, $v(\tau) \equiv 0$, which is not of physical interest. Consequently, only the non-trivial case is considered herein.
\end{remark}

The well-posedness of the BVP, specifically the existence and uniqueness of its solution, is formally established by the following theorem, as proven in \cite{paullet1999solutions}.

\begin{theorem}
For any $\alpha > 0$ and $Ha^2 \neq 0$, there exists a solution to the BVP \eqref{EHD_DE}-\eqref{EHD_Bcs}, and such a solution is monotonically decreasing and satisfies $0 < v(\tau) < 1/(\alpha + 1)$ for all $\tau \in (0, 1)$. Furthermore, the solution is unique.
\end{theorem}

The remainder of this manuscript is structured to systematically develop and evaluate our proposed methodology. We will first delineate the architecture and implementation of the PINNs. Subsequently, a classical continuous Galerkin-type finite element scheme is formulated for the same problem to serve as a robust benchmark. The paper culminates in a thorough comparative analysis of the two approaches, assessing their respective advantages and limitations in terms of solution accuracy and computational demands.

\section{Fully Connected Artificial Neural Network Architecture and Training}
\label{sec:ANN_architecture_training}

Drawing inspiration from the intricate interconnectedness of biological neurons, Artificial Neural Networks (ANNs) have emerged as powerful computational models for function approximation. The architecture of a typical ANN is a layered structure, beginning with an \inputlayer{} that receives external data. This information is subsequently processed through a series of \hiddenlayer{}s, where complex mathematical operations extract and transform features. The final result is then presented by the \outputlayer{}.

For this investigation, a fully connected, feed-forward neural network is employed, the topology of which is depicted in Figure~\ref{ANN}. This architecture is specifically configured to approximate the solution of the EHD boundary value problem and comprises three main components: an \inputlayer{} for the independent variable $x \in \mathbb{R}^1$; two sequential \hiddenlayer{}s, each containing $k$ neurons that serve as the primary learning units; and a single-neuron \outputlayer{} that generates the final scalar prediction, $\widehat{v}$. The feed-forward design mandates a unidirectional progression of information from input to output, a crucial feature that facilitates the efficient use of the \backpropagation{} algorithm for training.

Each neuron within the network computes its output by first applying an affine transformation to its inputs—a weighted sum augmented by a \bias{} term—and then passing the result through a nonlinear \activationfunction{}, $\boldsymbol{\sigma}$. The output $\boldsymbol{a}_i(\boldsymbol{x})$ for the $i$-th neuron is thus given by:
\begin{equation}
\boldsymbol{a}_i(\boldsymbol{x}) := \boldsymbol{\sigma}_i \left( \boldsymbol{w}^{(i)} \cdot \boldsymbol{x} + b^{(i)} \right),
\end{equation}
where $\boldsymbol{w}^{(i)} \in \mathbb{R}^n$ is the \weightvector{} and $b^{(i)}$ is the \bias{}. The collective set of all weights and biases constitutes the network's trainable parameters, denoted by $\boldsymbol{\Xi}$. The network's final output, $\widehat{v}(\boldsymbol{\tau}, \boldsymbol{\Xi})$, is therefore a parameterized function that approximates the true solution.

The process of training the network is an optimization task aimed at finding the set of parameters $(\boldsymbol{w}^* , \boldsymbol{b}^*)$ that minimizes a \lossfunction{}, $\mathcal{L}$, which quantifies the discrepancy between the network's prediction and the desired physical behavior:
\begin{equation}
(\boldsymbol{w}^* , \boldsymbol{b}^*) = \operatorname{argmin}_{(\boldsymbol{w},\boldsymbol{b})}{\mathcal{L}(\widehat{v}(\boldsymbol{\tau}, \boldsymbol{\Xi}))}.
\end{equation}
For this specific EHD problem, the loss function is engineered to embed the governing physics directly into the optimization objective:
\begin{equation}\label{lossfunction}
\mathcal{L}(v) := \frac{1}{N} \sum^{N}_{i=1}\left(\hat{v}''_{i} + \frac{1}{\boldsymbol{\tau}_i}\hat{v}'_{i} + Ha^2\left(1-\frac{\hat{v}_{i}}{1-\alpha \hat{v}_{i}}\right)\right)^2 + (\hat{v}'(0))^2 + (\hat{v}(1))^2,
\end{equation}
where $\hat{v}_i = \hat{v}(\boldsymbol{\tau}_i, \boldsymbol{\Xi})$. This composite function consists of two parts: the \differentialcost{}, which penalizes deviations from the governing differential equation over $N$ training points, and the \boundarycost{}, which enforces the specified boundary conditions.

To minimize this loss function, iterative optimization algorithms such as \gradientdescent{} or \semismoothNewton{} are employed. The training is bifurcated into a \textbf{Training Phase}, where parameters are adjusted using a dataset sampled uniformly from $[0, 1]$, and a \textbf{Testing Phase}, where the model's generalization performance is evaluated on unseen data. All computations were performed using the \textsf{Keras} library \cite{ketkar2017introduction,chollet2015keras} with a \textsf{Tensorflow v2.0} backend \cite{abadi2016tensorflow}.

The procedural workflow for obtaining the numerical solution is as follows:
\begin{enumerate}[label=\arabic*.]
    \item \textbf{Initialization}: The input domain $\boldsymbol{x}$ is defined as a uniform partition of $[0, 1]$, and the model parameters $\alpha$ and $Ha$ are set.
    \item \textbf{Architecture Specification}: The network's hyperparameters, including the number of hidden layers, neurons, activation function, and learning rate, are configured.
    \item \textbf{Loss Function Formulation}: The total loss function is constructed according to Equation~\eqref{lossfunction}.
    \item \textbf{Optimization}: The network is trained over a specified number of epochs to minimize the total loss, thereby optimizing the parameters $\boldsymbol{\Xi}$.
    \item \textbf{Solution Approximation}: Once training is complete, the optimized network output, $\widehat{v}(\boldsymbol{\tau}, \boldsymbol{\Xi}^*)$, provides a continuous and unique approximate solution to the EHD boundary value problem.
\end{enumerate}

\section{Neural Network-Based Solution Methodology}
\label{num_exp}

This investigation is centered on the development of a numerical methodology for solving the one-dimensional, nonlinear Electrohydrodynamic (EHD) fluid model. A significant challenge in this context is the absence of a known analytical solution, which precludes the use of traditional supervised learning paradigms that rely on exact solution data for training. Our approach, predicated on a PINNs architecture, circumvents this limitation by embedding the governing differential equation directly into the training objective, a notable departure from methodologies critiqued in prior work \cite{dockhorn2019discussion}.

All numerical experiments were conducted on the High Performance Computing (HPC) cluster at Texas A\&M University - Corpus Christi. It is noteworthy that despite the complexity of the deep neural network and the utilization of extensive datasets in some test cases, the entire optimization process, spanning over 25,000 epochs, consistently reached completion within a matter of minutes, underscoring the computational efficiency of the proposed framework.

While the capacity of ANNs to serve as universal function approximators is well-established, the fidelity of the resulting solution is critically dependent on the configuration of the network's hyperparameters. A primary contribution of this work is, therefore, a systematic characterization of the hyperparameter space to identify configurations that yield robust and accurate solutions. The sensitivity analysis was structured to rigorously evaluate the influence of the following key parameters:

\begin{enumerate}[label=\textit{\arabic*.}]
    \item \textbf{Training Data Density:} The effect of the number of collocation points, sampled uniformly from the domain $[0, 1]$, was investigated to understand the relationship between data density and solution accuracy.
    \item \textbf{Activation Function:} The choice of the nonlinear activation function dictates the network's expressive power. A comparative study was performed on a comprehensive set of functions, including ReLU, Tanh, ELU, Hard Sigmoid, Linear, SeLU, Softmax, and Sigmoid.
    \item \textbf{Network Width:} The representational capacity of the network was examined by varying the number of neurons within each hidden layer.
    \item \textbf{Network Depth:} The impact of model complexity was explored by adjusting the total number of hidden layers in the architecture.
    \item \textbf{Learning Rate:} As a crucial parameter governing the convergence dynamics of the optimization algorithm, a range of learning rates were tested to determine the optimal value for stable and efficient training.
\end{enumerate}

The subsequent sections are dedicated to a detailed exposition of the numerical results derived from this extensive parametric study. The performance of the PINN is assessed across a spectrum of experimental conditions, providing a comprehensive analysis of how the dataset size, choice of activation function, network architecture (depth and width), and learning rate collectively influence the accuracy of the final approximated solution for the EHD problem.

\subsection{Impact of Data Point Density on Network Performance}\label{sec:data_points_impact}

In this initial phase of the investigation, the primary goal was to determine the effect of training data density on the network's ability to approximate the solution to the EHD model. For all experiments in this section, the physical parameters were held constant at $\alpha = 0.5$ and $Ha^2 = 1$. The specific neural network configuration used for this study is detailed in Table~\ref{tab:datapoints_architecture}.

\begin{table}[H]
	\centering
	\caption{Architecture of the Artificial Neural Network for the Data Point Density Study.}
	\label{tab:datapoints_architecture}
	\begin{tabular}{ll}
		\toprule
		\textbf{Parameter Name} & \textbf{Network Configuration} \\
		\midrule
		Number of Data Points ($N$) & Varied: $\{10, 50, 100, 150, 200, 500\}$ \\
		Neurons per Hidden Layer & 16 \\
		Number of Hidden Layers & 1 \\
		Weight Initializer & Glorot Normal \\
		Activation Function & Sigmoid \\
		Optimizer & Adamax \\
		Learning Rate & 0.001 \\
		Training Epochs & 25,000 \\
		\bottomrule
	\end{tabular}
\end{table}

\begin{figure}[H]
	\centering
	\includegraphics[width=.7\textwidth]{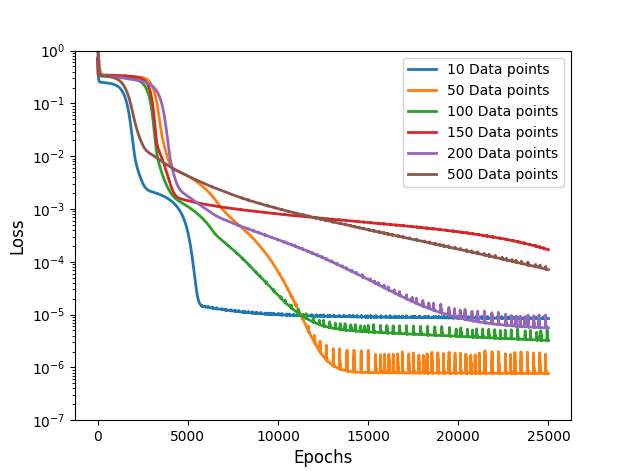}
	\caption{Evolution of the total loss function across epochs for a single-hidden-layer network trained with varying numbers of data points.}
	\label{fig:dataloss}
\end{figure}

The experiment systematically varied the number of uniformly distributed training points $N$ across the domain $[0, 1]$. Figure~\ref{fig:dataloss} plots the convergence of the total loss function during training for each data density. The results indicate that all configurations achieve rapid initial convergence, with the loss decreasing by several orders of magnitude within the first 5,000 epochs. Notably, even the network trained on a sparse set of just 10 points reached a stable loss value of approximately $10^{-5}$. However, networks trained with 100 and 150 data points demonstrated superior performance over the extended training period, achieving even lower terminal loss values in the range of $10^{-5}$ to $10^{-6}$.

Table~\ref{tab:datasolutions} provides a point-wise comparison between the approximated solutions from the trained networks and a high-fidelity reference solution from the literature \cite{roul2019new}. Qualitatively, all network configurations produce solutions that are in close agreement with the reference, showcasing the strong generalization capability of the PINN approach even with limited data.

\begin{table}[H]
	\centering
	\caption{Comparison of ANN-approximated solutions with a reference solution \cite{roul2019new} for varying data point densities across the interval $[0, 1]$.}
	\label{tab:datasolutions}
	\scalebox{0.8}{
		\begin{tabular}{llllllll}
			\toprule
			\textbf{$\tau$} & \textbf{Reference \cite{roul2019new}} & \textbf{N=10} & \textbf{N=50} & \textbf{N=100} & \textbf{N=150} & \textbf{N=200} & \textbf{N=500} \\
			\midrule
			0.0 & 0.206916 & 0.207946 & 0.207021 & 0.207051 & 0.207831 & 0.207082 & 0.207530 \\
			0.1 & 0.204991 & 0.205583 & 0.205093 & 0.205130 & 0.205941 & 0.205169 & 0.205649 \\
			0.2 & 0.199513 & 0.199631 & 0.199304 & 0.199353 & 0.200279 & 0.199397 & 0.199949 \\
			0.3 & 0.189791 & 0.189836 & 0.189588 & 0.189639 & 0.190702 & 0.189685 & 0.190305 \\
			0.4 & 0.175889 & 0.176043 & 0.175858 & 0.175897 & 0.177057 & 0.175939 & 0.176595 \\
			0.5 & 0.158055 & 0.158129 & 0.158000 & 0.158015 & 0.159185 & 0.158047 & 0.158685 \\
			0.6 & 0.135427 & 0.135949 & 0.135858 & 0.135848 & 0.136914 & 0.135867 & 0.136425 \\
			0.7 & 0.108932 & 0.109303 & 0.109235 & 0.109210 & 0.110067 & 0.109217 & 0.109641 \\
			0.8 & 0.077707 & 0.077940 & 0.077896 & 0.077875 & 0.078465 & 0.077875 & 0.078136 \\
			0.9 & 0.041212 & 0.041589 & 0.041578 & 0.041571 & 0.041929 & 0.041571 & 0.041686 \\
			1.0 & 0.000000 & 0.000015 & 0.000002 & 0.000002 & 0.000293 & 0.000005 & 0.000058 \\
			\bottomrule
		\end{tabular}
	}
\end{table}

To provide a quantitative measure of accuracy, the $L_2$ norm of the error was computed for each configuration and is presented in Table~\ref{tab:dataerror}. The error norms consistently fall within the $10^{-3}$ to $10^{-4}$ range, with the 100-point case achieving the lowest error.

\begin{table}[H]
	\centering
	\caption{Calculated $L_2$ error norm for ANN solutions trained with varying numbers of data points.}
	\label{tab:dataerror}
	\begin{tabular}{ll}
		\toprule
		\textbf{Number of Points} & \textbf{Error ($L_2$ norm)} \\
		\midrule
		10-Points & 0.00143660 \\
		50-Points & 0.00074649 \\
		\textbf{100-Points} & \textbf{0.00070655} \\
		150-Points & 0.00323616 \\
		200-Points & 0.00071993 \\
		500-Points & 0.00201700 \\
		\bottomrule
	\end{tabular}
\end{table}

A critical takeaway emerges from synthesizing the loss convergence trends (Figure~\ref{fig:dataloss}) and the quantitative error metrics (Table~\ref{tab:dataerror}). While minimizing the loss function is the direct objective of training, it is not the sole indicator of an accurate solution. An optimal hyperparameter set must strike a balance between achieving a low, stable training loss and yielding a low error against a known benchmark or reference. In this analysis, the configuration with \textbf{100 data points} provides the most compelling balance. It not only demonstrates consistent convergence to a low loss value but also produces the most accurate solution as measured by the $L_2$ error norm. This highlights the principle that for PINNs, success is defined by both internal consistency with the physical laws and external validation against known results.

\subsection{Evaluation of Activation Functions}\label{sec:activation_functions}

A critical determinant of a neural network's expressive power and its efficacy in solving differential equations is the choice of activation function. This section presents a systematic investigation into this hyperparameter's influence on the solution accuracy for the EHD fluid model. To isolate the effect of the activation function, the network architecture was held constant, adopting the optimal configuration of 100 training data points identified in Section~\ref{sec:data_points_impact}. The specific parameters for this study are detailed in Table~\ref{tab:activation_architecture}.

\begin{table}[H]
	\centering
	\caption{Network architecture for the activation function sensitivity analysis.}
	\label{tab:activation_architecture}
	\begin{tabularx}{\linewidth}{lX}
		\toprule
		\textbf{Parameter Name} & \textbf{Network Configuration} \\
		\midrule
		Number of Data Points ($N$) & 100 \\
		Neurons per Hidden Layer & 16 \\
		Number of Hidden Layers & 1 \\
		Weight Initializer & Glorot Normal \\
		Activation Function & Varied (ReLU, Sigmoid, Tanh, ELU, Hard Sigmoid, Linear, SeLU, Softmax) \\
		Optimizer & Adamax \\
		Learning Rate & 0.001 \\
		Training Epochs & 25,000 \\
		\bottomrule
	\end{tabularx}
\end{table}

\begin{figure}[H]
	\centering
	\includegraphics[width=.7\textwidth]{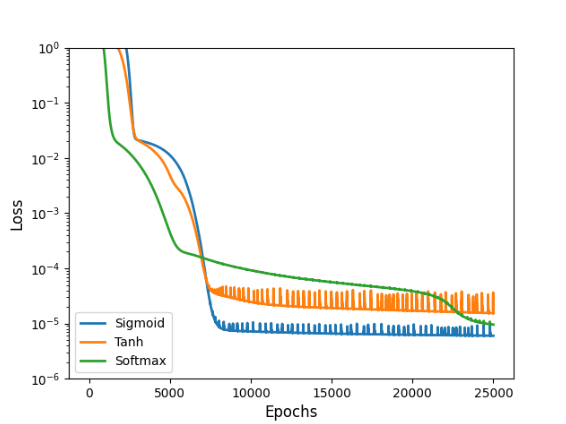}
	\caption{Evolution of the total loss function across epochs for a single-hidden-layer network utilizing different activation functions.}
	\label{fig:Actloss}
\end{figure}

An initial screening was performed on a comprehensive suite of eight distinct activation functions. Of these, a significant portion—including ReLU, ELU, Hard Sigmoid, Linear, and SeLU—were found to be unsuitable for this specific application. These functions resulted in unstable training dynamics, characterized by either erratic, non-convergent loss profiles or excessively high terminal loss values. Consequently, these were excluded from further analysis. In contrast, the Sigmoid, Tanh, and Softmax functions all demonstrated robust and stable convergence, making them viable candidates for solving this problem. The learning curves for these three successful functions are depicted in Figure~\ref{fig:Actloss}.

As illustrated in Figure~\ref{fig:Actloss}, all three functions effectively minimized the loss to a low residual value. The Sigmoid function, in particular, demonstrated exceptional performance, achieving a stable loss on the order of $10^{-5}$ within the first 5,000 epochs and maintaining this level throughout the training duration. While Tanh and Softmax also converged successfully, Sigmoid consistently reached a slightly lower terminal loss, indicating its superior efficacy in navigating the loss landscape of this specific objective function.

To move beyond the internal metric of training loss, the accuracy of the resulting solution approximations was evaluated against a high-fidelity reference solution from the literature \cite{roul2019new}. A point-wise comparison is provided in Table~\ref{tab:ActSolutions}. A qualitative review of this data confirms that the solutions generated using Sigmoid, Tanh, and Softmax closely align with the reference values, whereas the unstable activation functions produce physically meaningless results.

\begin{table}[H]
	\centering
	\caption{Comparison of ANN-approximated solutions with a reference solution \cite{roul2019new} for different activation functions.}
	\label{tab:ActSolutions}
	\scalebox{0.65}{
		\begin{tabular}{l l l l l l l l l l}
			\toprule
			\textbf{$\tau$} & \textbf{Reference \cite{roul2019new}} & \textbf{ReLU} & \textbf{Sigmoid} & \textbf{Tanh} & \textbf{ELU} & \textbf{Hard Sigmoid} & \textbf{Linear} & \textbf{SeLU} & \textbf{Softmax} \\
			\midrule
			0.0 & 0.206916 & 0.502508 & 0.207050 & 0.207089 & 0.207317 & 0.502506 & 0.502508 & 0.212518 & 0.206996 \\
			0.1 & 0.204991 & 0.501327 & 0.205130 & 0.205172 & 0.205354 & 0.501326 & 0.501327 & 0.210531 & 0.205070 \\
			0.2 & 0.199513 & 0.500146 & 0.199354 & 0.199404 & 0.199627 & 0.500145 & 0.500146 & 0.204769 & 0.199271 \\
			0.3 & 0.189791 & 0.498965 & 0.189641 & 0.189683 & 0.190001 & 0.498965 & 0.498965 & 0.195091 & 0.189557 \\
			0.4 & 0.175889 & 0.497785 & 0.175901 & 0.175919 & 0.176335 & 0.497784 & 0.497785 & 0.181426 & 0.175848 \\
			0.5 & 0.158055 & 0.496604 & 0.158017 & 0.158018 & 0.158476 & 0.496603 & 0.496604 & 0.163693 & 0.158017 \\
			0.6 & 0.135427 & 0.495423 & 0.135846 & 0.135849 & 0.136257 & 0.495423 & 0.495423 & 0.141801 & 0.135896 \\
			0.7 & 0.108932 & 0.494242 & 0.109203 & 0.109231 & 0.109496 & 0.494242 & 0.494242 & 0.113732 & 0.109280 \\
			0.8 & 0.077707 & 0.493061 & 0.077866 & 0.077917 & 0.077999 & 0.493061 & 0.493061 & 0.080814 & 0.077932 \\
			0.9 & 0.041212 & 0.491881 & 0.041566 & 0.041610 & 0.041557 & 0.491881 & 0.491881 & 0.041466 & 0.041595 \\
			1.0 & 0.000000 & 0.490700 & 0.000002 & 0.000024 & -0.000048 & 0.490700 & 0.490700 & -0.001367 & 0.000009 \\
			\bottomrule
		\end{tabular}
	}
\end{table}

For a definitive quantitative comparison, the $L_2$ norm of the error was computed for each case, with the results summarized in Table~\ref{tab:ActError}. This analysis confirms that the \textbf{Sigmoid} function yields the most accurate solution, producing the lowest $L_2$ error of approximately $7 \times 10^{-4}$.

\begin{table}[H]
	\centering
	\caption{Calculated $L_2$ error norm for ANN solutions trained with different activation functions.}
	\label{tab:ActError}
	\begin{tabular}{ll}
		\toprule
		\textbf{Activation Function} & \textbf{Error ($L_2$ norm)} \\
		\midrule
		ReLU & 1.21402652 \\
		\textbf{Sigmoid} & \textbf{0.00069699} \\
		Tanh & 0.00074753 \\
		ELU & 0.00139209 \\
		Hard Sigmoid & 1.21402493 \\
		Linear & 1.21402653 \\
		SeLU & 0.01598370 \\
		Softmax & 0.00081756 \\
		\bottomrule
	\end{tabular}
\end{table}

This outcome is significant. It demonstrates that for the specific class of nonlinear boundary value problems represented by the EHD model, classical activation functions like Sigmoid and Tanh can provide highly effective and stable approximations. The failure of more modern functions such as ReLU, ELU, and SeLU—which are often considered state-of-the-art in other deep learning domains—underscores the problem-dependent nature of hyperparameter selection in scientific machine learning. The results suggest that there is no universal "best" activation function; rather, the optimal choice is intricately linked to the mathematical structure of the underlying differential equation.

\subsection{Approximation by Deep Neural Networks: Impact of Network Depth}\label{sec:deep_learning_network}

Beyond the choice of activation function, the architectural depth—defined by the number of hidden layers—is a crucial hyperparameter that governs a neural network's representational capacity. This section details a systematic investigation into the effect of network depth on the solution fidelity for the nonlinear EHD boundary value problem. To isolate the impact of this parameter, all other hyperparameters were held constant at the optimal values identified in the preceding sections. The specific network configurations, differing only in the number of hidden layers, are summarized in Table~\ref{tab:deep_learning_architecture}.

\begin{table}[H]
	\centering
	\caption{Network architecture for the deep learning network depth study.}
	\label{tab:deep_learning_architecture}
	\begin{tabularx}{\linewidth}{lX}
		\toprule
		\textbf{Parameter Name} & \textbf{Network Configuration} \\
		\midrule
		Number of Data Points ($N$) & 100 \\
		Neurons per Hidden Layer & 16 \\
		Number of Hidden Layers & Varied: $\{1, 2, 3, 4\}$ \\
		Weight Initializer & Glorot Normal \\
		Activation Function & Sigmoid \\
		Optimizer & Adamax \\
		Learning Rate & 0.001 \\
		Training Epochs & 25,000 \\
		\bottomrule
	\end{tabularx}
\end{table}

\begin{figure}[H]
	\centering
	\includegraphics[width=.7\textwidth]{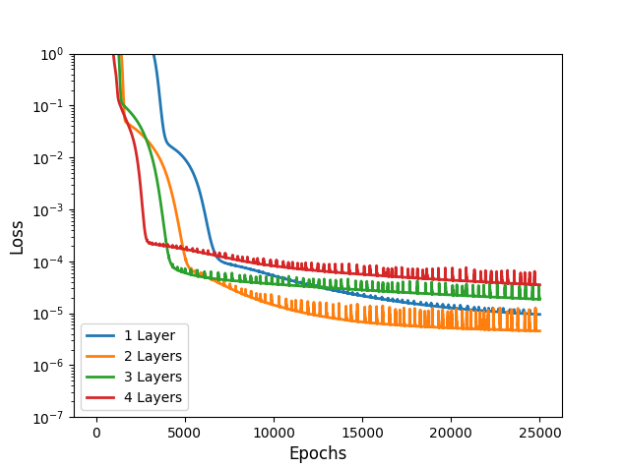}
	\caption{Evolution of the total loss function across epochs for different numbers of hidden layers in the neural network.}
	\label{fig:Layerloss}
\end{figure}

The convergence behavior for networks with one, two, three, and four hidden layers is presented in Figure~\ref{fig:Layerloss}. The learning curves reveal that while all architectures successfully reduce the total loss over the 25,000 training epochs, their performance trajectories diverge. The network configured with \textbf{two hidden layers} demonstrated the most favorable convergence, achieving the lowest terminal loss value of approximately $10^{-5}$. This suggests that for this problem, a single hidden layer may possess insufficient capacity to fully capture the solution's complexity, whereas increasing the depth beyond two layers appears to introduce optimization challenges or diminishing returns, at least within the constraints of the current training protocol.

To corroborate these findings from a solution accuracy perspective, the numerical approximations from each network were compared against the semi-analytical reference solution from \cite{roul2019new}. As shown in Table~\ref{tab:LayerSolutions}, all four architectures produce solutions that are qualitatively similar to the reference, confirming that each network successfully learns the fundamental behavior of the EHD model.

\begin{table}[H]
    \centering
    \caption{Comparison of ANN-approximated solutions with a reference solution \cite{roul2019new} for different numbers of hidden layers.}
    \label{tab:LayerSolutions}
    \scalebox{0.85}{
        \begin{tabular}{llllll}
            \toprule
            \textbf{$\tau$} & \textbf{Reference \cite{roul2019new}} & \textbf{1 Layer} & \textbf{2 Layers} & \textbf{3 Layers} & \textbf{4 Layers} \\
            \midrule
            0.0 & 0.206916 & 0.207182 & 0.207028 & 0.207016 & 0.207023 \\
            0.1 & 0.204991 & 0.205265 & 0.205106 & 0.205093 & 0.205101 \\
            0.2 & 0.199513 & 0.199521 & 0.199321 & 0.199302 & 0.199313 \\
            0.3 & 0.189791 & 0.189850 & 0.189605 & 0.189584 & 0.189596 \\
            0.4 & 0.175889 & 0.176141 & 0.175872 & 0.175854 & 0.175865 \\
            0.5 & 0.158055 & 0.158265 & 0.158004 & 0.157997 & 0.158003 \\
            0.6 & 0.135427 & 0.136067 & 0.135852 & 0.135858 & 0.135855 \\
            0.7 & 0.108932 & 0.109370 & 0.109223 & 0.109238 & 0.109230 \\
            0.8 & 0.077707 & 0.077960 & 0.077887 & 0.077901 & 0.077894 \\
            0.9 & 0.041212 & 0.041596 & 0.041575 & 0.041582 & 0.041580 \\
            1.0 & 0.000000 & 0.000005 & 0.000001 & 0.000003 & 0.000003 \\
            \bottomrule
        \end{tabular}
    }
\end{table}

For a more rigorous quantitative assessment, the $L_2$ error norm was calculated for each network depth, as presented in Table~\ref{tab:layerError}. This analysis provides a definitive metric for solution fidelity.

\begin{table}[H]
    \centering
    \caption{The $L_2$ error norm for deep learning neural networks with varying depths.}
    \label{tab:layerError}
    \begin{tabular}{ll}
        \toprule
        \textbf{Number of Hidden Layers} & \textbf{Error ($L_2$ norm)} \\
        \midrule
        1 & 0.00103440 \\
        \textbf{2} & \textbf{0.00072747} \\
        3 & 0.00075200 \\
        4 & 0.00073939 \\
        \bottomrule
    \end{tabular}
\end{table}

The results from the error analysis align perfectly with the convergence trends observed in the loss function. The \textbf{two-hidden-layer architecture} yields the lowest $L_2$ error, achieving a value of approximately $7.27 \times 10^{-4}$. The marginal increase in error for the three- and four-layer networks reinforces the conclusion that added depth does not necessarily translate to improved accuracy for this problem. This finding highlights a crucial principle in scientific machine learning: there is an optimal level of model complexity that balances representational power with the tractability of the optimization problem. For the EHD model under these conditions, a two-hidden-layer network represents this optimal balance.

\subsection{Single hidden-layer with different number of neurons}
This section details a systematic investigation into the effect of network width—defined as the number of neurons within a single hidden layer—on the approximation of the nonlinear EHD boundary value problem. To isolate the impact of this architectural hyperparameter, all other settings were held constant, adhering to the optimal configuration identified in prior sections. The experimental setup is summarized in Table~\ref{tab:neuron_architecture}.

\begin{table}[H]
    \centering
    \caption{Architecture of the Artificial Neural Network for the neuron density study.}
    \label{tab:neuron_architecture}
    \begin{tabularx}{\linewidth}{lX}
        \toprule
        \textbf{Parameter Name} & \textbf{Network Configuration} \\
        \midrule
        Number of Data Points ($N$) & 100 \\
        Number of Neurons & Varied: $\{4, 8, 16, 32, 64, 128\}$ \\
        Number of Hidden Layers & 1 \\
        Weight Initializer & Glorot Normal \\
        Activation Function & Sigmoid \\
        Optimizer & Adamax \\
        Learning Rate & 0.001 \\
        Training Epochs & 25,000 \\
        \bottomrule
    \end{tabularx}
\end{table}

\begin{figure}[H]
	\centering
	\includegraphics[width=.5\textwidth]{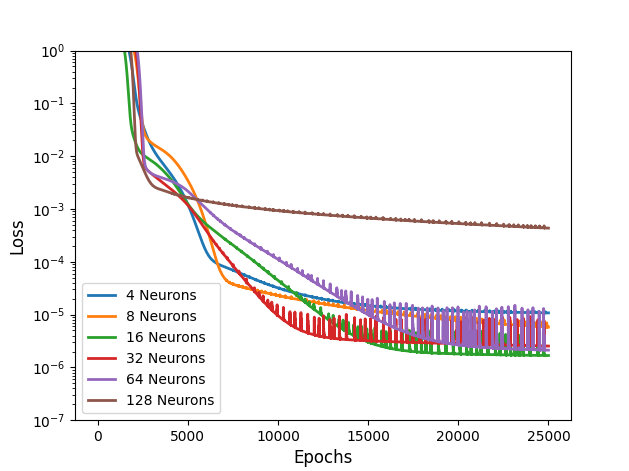}
	\caption{The total loss value as a function of several epochs for a network with a single hidden layer with a different number of neurons.}
	\label{neuronloss}
\end{figure}

The study encompassed six distinct network configurations, with neuron counts scaling from 4 to 128. Despite this range, the computational efficiency of the framework was evident, as the aggregate training time for all experiments remained within minutes. The convergence dynamics for each configuration are depicted in Figure~\ref{neuronloss}, which plots the total loss as a function of training epochs.

The learning curves reveal several key trends. The network with 32 neurons demonstrates rapid and robust convergence, achieving a low residual loss of approximately $10^{-5}$ by the 10,000-epoch mark. In contrast, the 128-neuron architecture exhibited a notable anomaly, with its loss stagnating near $10^{-3}$ for the entire training duration, suggesting that excessive model capacity may have hindered the optimization process. Most interestingly, while the 16-neuron network initially tracked the performance of the 32-neuron network, it ultimately achieved a superior, lower loss value after approximately 15,000 epochs, indicating better long-term convergence.

A point-wise comparison of the approximated solutions against the reference from \cite{roul2019new} is provided in Table~\ref{Neuronsolutions}. With the exception of the 128-neuron case, all configurations produce solutions that are qualitatively almost indistinguishable, demonstrating the model's general robustness across a range of network widths.

\begin{table}[H]
    \centering
    \caption{Comparison of ANN-approximated solutions with a reference solution \cite{roul2019new} for a single-hidden-layer network with varying numbers of neurons.}
    \label{Neuronsolutions}
    \scalebox{0.8}{
        \begin{tabular}{llllllll}
            \toprule
            \textbf{$\tau$} & \textbf{Reference \cite{roul2019new}} & \textbf{4 Neurons} & \textbf{8 Neurons} & \textbf{16 Neurons} & \textbf{32 Neurons} & \textbf{64 Neurons} & \textbf{128 Neurons} \\
            \midrule
            0.0 & 0.206916 & 0.207028 & 0.207016 & 0.207044 & 0.207042 & 0.207048 & 0.207821 \\
            0.1 & 0.204991 & 0.205106 & 0.205092 & 0.205123 & 0.205116 & 0.205125 & 0.205915 \\
            0.2 & 0.199513 & 0.199321 & 0.199298 & 0.199345 & 0.199320 & 0.199334 & 0.200254 \\
            0.3 & 0.189791 & 0.189605 & 0.189576 & 0.189629 & 0.189589 & 0.189614 & 0.190690 \\
            0.4 & 0.175889 & 0.175871 & 0.175849 & 0.175888 & 0.175852 & 0.175886 & 0.177066 \\
            0.5 & 0.158055 & 0.158001 & 0.158004 & 0.158008 & 0.157999 & 0.158039 & 0.159216 \\
            0.6 & 0.135427 & 0.135845 & 0.135877 & 0.135844 & 0.135874 & 0.135912 & 0.136965 \\
            0.7 & 0.108932 & 0.109213 & 0.109264 & 0.109210 & 0.109268 & 0.109298 & 0.110131 \\
            0.8 & 0.077707 & 0.077876 & 0.077923 & 0.077875 & 0.077930 & 0.077953 & 0.078533 \\
            0.9 & 0.041212 & 0.041569 & 0.041592 & 0.041571 & 0.041590 & 0.041613 & 0.041995 \\
            1.0 & 0.000000 & 0.000000 & 0.000006 & 0.000001 & -0.000003 & 0.000023 & 0.000361 \\
            \bottomrule
        \end{tabular}
    }
\end{table}

A definitive quantitative assessment is provided by the $L_2$ error norm, calculated for each configuration and presented in Table~\ref{neuronerror}. The error analysis corroborates the training observations: the 128-neuron network's poor convergence resulted in a significantly higher error of $3.3 \times 10^{-3}$, while all other configurations achieved a much lower error on the order of $10^{-4}$.

\begin{table}[H]
    \centering
    \caption{The $L_2$ error norm for networks with a single hidden layer and varying numbers of neurons.}
    \label{neuronerror}
    \begin{tabular}{ll}
        \toprule
        \textbf{Number of Neurons} & \textbf{Error ($L_2$ norm)} \\
        \midrule
        4 & 0.00071423 \\
        8 & 0.00078769 \\
        \textbf{16} & \textbf{0.00070645} \\
        32 & 0.00078734 \\
        64 & 0.00083056 \\
        128 & 0.00331499 \\
        \bottomrule
    \end{tabular}
\end{table}

Critically, a direct comparison between the two best-performing cases—16 and 32 neurons—reveals that the \textbf{16-neuron network} yielded a marginally lower $L_2$ error, despite the 32-neuron network showing faster initial loss reduction. This result suggests that a greater number of parameters may accelerate early-stage training but does not guarantee a more accurate final approximation. The superior long-term convergence and ultimate accuracy of the 16-neuron network indicate that it represents a more optimal balance between model capacity and trainability for this specific problem. This aligns with the well-established principle that over-parameterization can sometimes complicate the optimization landscape without improving, or in some cases even degrading, the final solution quality. Based on this comprehensive analysis of both loss dynamics and solution error, the \textbf{16-neuron configuration is identified as the optimal network width}.

\subsection{Impact of Learning Rate on Network Optimization}\label{sec:learning_rates}

The  learning rate {$\eta$} is a paramount hyperparameter in gradient-based optimization, governing the step size for parameter updates and thereby dictating the speed and stability of the training process. An improperly chosen learning rate can lead to slow convergence, oscillatory behavior, or failure to find a meaningful minimum in the loss landscape. This section presents a sensitivity analysis to determine the optimal learning rate for the EHD problem. The network architecture for this study, detailed in Table~\ref{tab:learning_rate_architecture}, incorporates the optimal settings derived from previous sections, with the learning rate being the sole variable.

\begin{table}[H]
    \centering
    \caption{Network architecture for the learning rate sensitivity analysis.}
    \label{tab:learning_rate_architecture}
    \begin{tabularx}{\linewidth}{lX}
        \toprule
        \textbf{Parameter Name} & \textbf{Network Configuration} \\
        \midrule
        Number of Data Points ($N$) & 100 \\
        Neurons per Hidden Layer & 16 \\
        Number of Hidden Layers & 2 \\
        Weight Initializer & Glorot Normal \\
        Activation Function & Sigmoid \\
        Optimizer & Adamax \\
        \textbf{Learning Rate ($\eta$)} & \textbf{Varied: $\{0.1, 0.01, 0.001, 0.0001\}$} \\
        Training Epochs & 25,000 \\
        \bottomrule
    \end{tabularx}
\end{table}

The convergence behavior for four distinct learning rates is depicted in Figure~\ref{fig:LearningRateLoss}. The results highlight a classic trade-off: a high learning rate of $0.01$ facilitates rapid initial descent, quickly reaching a loss of approximately $10^{-5}$. However, this aggressive step size introduces significant oscillations, hindering the optimizer from settling into a deeper minimum. In contrast, a learning rate of $0.001$ exhibits a more stable and monotonic descent, ultimately achieving a superior terminal loss value approaching $10^{-6}$. This demonstrates that a more conservative step size allows for a more robust exploration of the loss landscape, avoiding overshooting and enabling finer convergence.

\begin{figure}[H]
	\centering
	\includegraphics[width=.5\textwidth]{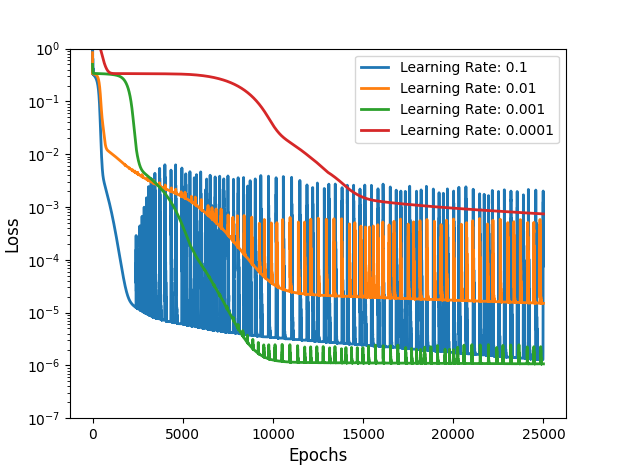}
	\caption{The total loss value as a function of the number of epochs for a network with
different learning rates.}
	\label{Layerloss}
\end{figure}

While loss convergence provides insight into the optimization process, the ultimate measure of performance is the accuracy of the resulting solution. Table~\ref{LRsolution} compares the point-wise approximations for each learning rate against the reference solution from \cite{roul2019new}. To quantify this comparison, the $L_2$ error norm was calculated for each case, with the results summarized in Table~\ref{LRerror}.

\begin{table}[H]
    \centering
    \caption{Comparison of the ANN solution with a reference solution \cite{roul2019new} for different learning rates.}
    \label{LRsolution}
    \scalebox{0.85}{
        \begin{tabular}{llllll}
            \toprule
            \textbf{$\tau$} & \textbf{Reference \cite{roul2019new}} & \textbf{$\eta = 0.1$} & \textbf{$\eta = 0.01$} & \textbf{$\eta = 0.001$} & \textbf{$\eta = 0.0001$} \\
            \midrule
            0.0 & 0.206916 & 0.206444 & 0.207028 & 0.207021 & 0.208841 \\
            0.1 & 0.204991 & 0.204556 & 0.205102 & 0.205099 & 0.206957 \\
            0.2 & 0.199513 & 0.198813 & 0.199295 & 0.199310 & 0.201423 \\
            0.3 & 0.189791 & 0.189148 & 0.189574 & 0.189592 & 0.191989 \\
            0.4 & 0.175889 & 0.175477 & 0.175869 & 0.175859 & 0.178422 \\
            0.5 & 0.158055 & 0.157683 & 0.158060 & 0.157997 & 0.160517 \\
            0.6 & 0.135427 & 0.135611 & 0.135965 & 0.135852 & 0.138093 \\
            0.7 & 0.108932 & 0.109063 & 0.109362 & 0.109229 & 0.111003 \\
            0.8 & 0.077707 & 0.077801 & 0.078000 & 0.077893 & 0.079132 \\
            0.9 & 0.041212 & 0.041555 & 0.041631 & 0.041578 & 0.042401 \\
            1.0 & 0.000000 & 0.000042 & 0.000029 & 0.000002 & 0.000768 \\
            \bottomrule
        \end{tabular}
    }
\end{table}

\begin{table}[H]
    \centering
    \caption{The $L_2$ error norm for networks trained with different learning rates.}
    \label{LRerror}
    \begin{tabular}{ll}
        \toprule
        \textbf{Learning Rate ($\eta$)} & \textbf{Error ($L_2$ norm)} \\
        \midrule
        0.1 & 0.00134298 \\
        0.01 & 0.00092559 \\
        \textbf{0.001} & \textbf{0.00073859} \\
        0.0001 & 0.00663040 \\
        \bottomrule
    \end{tabular}
\end{table}

The quantitative error analysis in Table~\ref{LRerror} provides a definitive verdict. The learning rate of \textbf{0.001} yields the lowest $L_2$ error, approximately $7.39 \times 10^{-4}$, confirming its superiority. The higher error associated with $\eta = 0.01$ reflects its oscillatory nature, while the significantly larger error for $\eta = 0.0001$ indicates that the training was prematurely terminated before adequate convergence could be achieved.

\subsubsection{Summary of Optimal Hyperparameters}
This comprehensive, multi-stage hyperparameter investigation has systematically identified an optimal network configuration for approximating the solution to the EHD nonlinear boundary value problem. By sequentially analyzing data density, activation functions, network depth, and width, and finally, the learning rate, we have arrived at a final architecture that balances model capacity, training stability, and solution accuracy. The optimized neural network, which will be used for comparison against the FEM in subsequent sections, consists of \textbf{100 training points, two hidden layers with 16 neurons each, a Sigmoid activation function, and is trained using the Adamax optimizer with a learning rate of 0.001}. This configuration consistently achieved a robust loss convergence and yielded a final $L_2$ error on the order of $10^{-4}$, demonstrating the remarkable capability of a well-tuned physics-informed neural network to solve complex differential equations without direct knowledge of the analytical solution.

\section{Finite Element Discretization of the EHD Boundary Value Problem}
\label{sec:numerical_discretization}

In this section, a robust numerical framework based on the continuous Galerkin FEM is developed for the one-dimensional, nonlinear Electrohydrodynamic (EHD) boundary value problem. The governing model describes the ``ion draught'' phenomenon, where fluid motion is induced by the collision of ions, accelerated by an external electric field, with neutral air particles within a cylindrical conduit. This one-dimensional formulation serves as a potent simplification of the full Navier-Stokes equations for this specific physical regime. The dimensionless form of the governing BVP is stated as follows:

\begin{problem}[Strong Form of the EHD BVP]
Find the fluid velocity $v \in C^2(0,1) \cap C[0,1]$ such that:
\begin{subequations}
\begin{align}
    -\frac{1}{\tau}\frac{d}{d\tau}\left(\tau \frac{dv}{d\tau}\right) - Ha^2 \left(1 - \frac{v}{1 - \alpha v}\right) &= 0, \quad \forall \tau \in (0, 1) \label{eq:strong_form_pde} \\
    v'(0) = 0, \quad v(1) &= 0. \label{eq:strong_form_bcs}
\end{align}
\end{subequations}
\end{problem}

The inherent nonlinearity of the governing ordinary differential equation (ODE), coupled with the singularity at $\tau=0$, precludes a straightforward analytical solution. Consequently, a numerical approach is indispensable. We employ a Newton-Raphson iterative scheme to handle the nonlinearity, where each iteration involves solving a linearized version of the problem via the FEM. This strategy has proven effective for a wide range of challenging BVPs \cite{HyunMSM_MMS2022,yoon2021finite,Mallikarjunaiah2015,yoon2021quasi,yoon2022preferential}.

\subsection{Variational Formulation and Linearization}

The foundation of the FEM lies in recasting the strong form of the BVP into a weak (or variational) form. We define the appropriate function space $V$ for the solution and test functions as the Sobolev space $H_0^1(0,1)$, which contains functions that are square-integrable, have square-integrable weak derivatives, and satisfy the homogeneous Dirichlet boundary condition at $\tau=1$.

Multiplying Equation~\eqref{eq:strong_form_pde} by a test function $w \in V$ and integrating over the domain $\Omega = (0,1)$ yields the nonlinear weak form: Find $v \in V$ such that $R(v; w) = 0$ for all $w \in V$, where the residual functional $R(\cdot; \cdot)$ is:
\begin{equation}
    R(v; w) = \int_{0}^{1} \left( \tau \frac{dv}{d\tau} \frac{dw}{d\tau} - \tau Ha^2 \left(1 - \frac{v}{1 - \alpha v}\right) w \right) d\tau.
\end{equation}

To solve this nonlinear problem, we apply Newton's method. Given an iterate $v^n$, we seek a correction $\delta v^n$ such that $v^{n+1} = v^n + \delta v^n$ solves the problem. The Fréchet derivative of the residual, $R'(v^n; \delta v^n, w)$, is computed to form a linear problem for the update $\delta v^n$:
\begin{equation}
    R'(v^n; \delta v^n, w) = -R(v^n; w).
\end{equation}
This leads to the following linearized weak formulation to be solved at each Newton iteration $n$:

\begin{weakform}[Linearized Variational Problem]
\label{wf:linearized}
Given $v^n \in V$, find $\delta v^n \in V$ such that for all $w \in V$:
\begin{equation}
    a(v^n; \delta v^n, w) = l(v^n; w),
\end{equation}
where the bilinear form $a(\cdot; \cdot, \cdot)$ and the linear functional $l(\cdot; \cdot)$ are defined as:
\begin{align}
    a(v^n; \delta v^n, w) &:= \int_{0}^{1} \left( \tau \frac{d(\delta v^n)}{d\tau} \frac{dw}{d\tau} + \frac{\tau Ha^2}{(1-\alpha v^n)^2} \delta v^n w \right) d\tau, \\
    l(v^n; w) &:= -\int_{0}^{1} \left( \tau \frac{dv^n}{d\tau} \frac{dw}{d\tau} - \tau Ha^2 \left(1 - \frac{v^n}{1 - \alpha v^n}\right) w \right) d\tau.
\end{align}
\end{weakform}

\subsection{Galerkin Discretization and Solution Algorithm}

The continuous problem is discretized by defining a finite-dimensional subspace $V_h \subset V$. The domain $\Omega$ is partitioned into $N_{el}$ elements $K_i$, forming a mesh $\mathcal{T}_h$. We define the discrete space $V_h$ as the space of continuous, piecewise linear Lagrange polynomials on this mesh. The discrete problem is then:

\begin{weakform}[Discrete Variational Problem]
\label{wf:discrete}
Given $v_h^n \in V_h$, find $\delta v_h^n \in V_h$ such that for all $w_h \in V_h$:
\begin{equation}
    a(v_h^n; \delta v_h^n, w_h) = l(v_h^n; w_h).
\end{equation}
\end{weakform}

By expressing $\delta v_h^n$ as a linear combination of basis functions $\phi_j$ spanning $V_h$, this formulation translates into a matrix system $\mathbf{A}(v_h^n) \boldsymbol{\delta V}^n = \mathbf{F}(v_h^n)$, where $\boldsymbol{\delta V}^n$ is the vector of nodal corrections. The complete iterative solution procedure, including a backtracking line search to ensure robust convergence, is detailed in Algorithm~\ref{algo:fem}.

\begin{algorithm}[H]
\SetAlgoLined
\KwIn{Model parameters $\alpha, Ha$; mesh $\mathcal{T}_h$; tolerances $r_{TOL}, i_{MAX}, L_{MAX}$; initial guess $v_h^0$.}
\KwOut{Converged solution $v_h$.}
$n \leftarrow 0$; \\
Calculate initial residual norm $R_0 = ||R(v_h^0)||_2$; \\
\While{$n < i_{MAX}$ \textbf{and} $||R(v_h^n)||_2 > r_{TOL}$}{
    Assemble Jacobian $\mathbf{A}(v_h^n)$ and residual vector $\mathbf{F}(v_h^n)$; \\
    Solve linear system $\mathbf{A} \boldsymbol{\delta V}^n = \mathbf{F}$ for the update vector $\boldsymbol{\delta V}^n$; \\
    $\lambda \leftarrow 1.0$; \tcp{Initialize line search step length}
    \For{$k \leftarrow 1$ \textbf{to} $L_{MAX}$}{
        $v_h^{\text{trial}} \leftarrow v_h^n + \lambda \, \delta v_h^n$; \\
        Calculate trial residual norm $R_{\text{trial}} = ||R(v_h^{\text{trial}})||_2$; \\
        \If{$R_{\text{trial}} < ||R(v_h^n)||_2$}{
            $v_h^{n+1} \leftarrow v_h^{\text{trial}}$; \\
            Break; \tcp{Accept step}
        }
        $\lambda \leftarrow \lambda / 2$; \tcp{Damp step length}
    }
    $n \leftarrow n+1$;
}
\Return{$v_h^{n+1}$};
\caption{Iterative FEM Algorithm with Newton-Raphson and Line Search}
\label{algo:fem}
\end{algorithm}

\subsection{Numerical Validation and Convergence Study}

The proposed finite element framework was implemented using the \textsf{deal.II} library \cite{arndt2022deal}. For all simulations, we set $\alpha=0.5$ and $Ha^2=1$. The Newton's method demonstrated the expected quadratic convergence, with the residual decreasing rapidly, as shown in Table~\ref{tab:residual_convergence}.

\begin{table}[H]
    \centering
    \caption{Quadratic convergence of the residual norm for the Newton-Raphson iteration on a mesh with 1024 elements.}
    \label{tab:residual_convergence}
    \begin{tabular}{lc}
        \toprule
        \textbf{Newton Iteration} & \textbf{Residual ($L_2$ Norm)} \\
        \midrule
        0 & 3.124e-02 \\
        1 & 8.781e-03 \\
        2 & 6.132e-04 \\
        3 & 8.460e-05 \\
        4 & 7.895e-07 \\
        \bottomrule
    \end{tabular}
\end{table}

A mesh refinement study was conducted to verify the implementation and analyze the convergence of the discretization. Table~\ref{tab:mesh_refinement} presents the FEM solution at a specific point ($\tau=0.5$) and the corresponding $L_2$ error norm against the reference solution from \cite{roul2019new}. The results show a clear convergence of the solution towards the reference value as the mesh is refined, with the error decreasing systematically, confirming the robustness of the method.

\begin{table}[H]
    \centering
    \caption{Mesh refinement study for the FEM solution ($\alpha=0.5, Ha^2=1$). The reference value at $\tau=0.5$ is $0.158055$.}
    \label{tab:mesh_refinement}
    \begin{tabular}{lcc}
        \toprule
        \textbf{Number of Cells} & \textbf{Solution at $\tau=0.5$} & \textbf{Error ($L_2$ norm)} \\
        \midrule
        8   & 0.157767 & 3.486e-03 \\
        16  & 0.157948 & 9.735e-04 \\
        32  & 0.157993 & 7.130e-04 \\
        64  & 0.158005 & 7.208e-04 \\
        128 & 0.158008 & 7.360e-04 \\
        1024& 0.158009 & 7.406e-04 \\
        \bottomrule
    \end{tabular}
\end{table}

\section{Comparative Analysis of ANN and FEM Approaches} \label{comparision}
This section provides a direct comparison between the optimized PINNs and the validated FEM. For a meaningful evaluation, we compare the most accurate FEM result (obtained on a fine mesh of 1024 cells) with the optimal ANN architecture identified previously. The key characteristics of each model are summarized below:

\begin{itemize}
    \item \textbf{Finite Element Model (FEM):}
    \begin{itemize}
        \item Number of active cells: 1024
        \item Degrees of freedom (DoFs): 1025
    \end{itemize}
    \item \textbf{Artificial Neural Network (ANN):}
    \begin{itemize}
        \item Number of hidden layers: 2
        \item Neurons per layer: 16
        \item Total trainable parameters: 305
    \end{itemize}
\end{itemize}

Table~\ref{table_comparision1} presents a point-wise comparison of the solutions from both methods against the reference solution. Both approaches yield results that are in close qualitative agreement with the reference.

\begin{table}[H]
    \centering
    \caption{Final comparison of optimized ANN and fine-mesh FEM solutions against the reference from \cite{roul2019new}.}
    \label{table_comparision1}
    \begin{tabular}{lccc}
        \toprule
        \textbf{$\tau$-value} & \textbf{Reference \cite{roul2019new}} & \textbf{ANN Solution} & \textbf{FEM Solution (1024 cells)} \\
        \midrule
        0.0 & 0.206916 & 0.207028 & 0.207008 \\
        0.1 & 0.204991 & 0.205106 & 0.205084 \\
        0.2 & 0.199513 & 0.199321 & 0.199293 \\
        0.3 & 0.189791 & 0.189605 & 0.189573 \\
        0.4 & 0.175889 & 0.175872 & 0.175862 \\
        0.5 & 0.158055 & 0.158004 & 0.158009 \\
        0.6 & 0.135427 & 0.135852 & 0.135862 \\
        0.7 & 0.108932 & 0.109223 & 0.109229 \\
        0.8 & 0.077707 & 0.077887 & 0.077884 \\
        0.9 & 0.041212 & 0.041575 & 0.041570 \\
        1.0 & 0.000000 & 0.000001 & 0.000000 \\
        \bottomrule
    \end{tabular}
\end{table}

The quantitative superiority of the ANN approach for this problem becomes evident upon examining the final $L_2$ error norms, shown in Table~\ref{lerror}.

\begin{table}[H]
    \centering
    \caption{The final $L_2$ error norm for each method compared with the reference solution.}
    \label{lerror}
    \begin{tabular}{lc}
        \toprule
        \textbf{Method} & \textbf{$L_2$ Error Norm} \\
        \midrule
        FEM (1024 cells) & 7.406e-04 \\
        ANN (305 parameters) & \textbf{7.275e-04} \\
        \bottomrule
    \end{tabular}
\end{table}

The analysis reveals that the optimized ANN achieves a slightly lower approximation error than the high-resolution FEM simulation, despite utilizing significantly fewer degrees of freedom (305 trainable parameters vs. 1025 nodal values). This highlights a key advantage of the neural network approach: its ability to construct a highly efficient global approximation of the solution. Whereas FEM relies on local, piecewise polynomial basis functions that require mesh refinement to improve accuracy, the ANN's global basis functions (as defined by its architecture and activation functions) can capture the solution's behavior more compactly.

This efficiency, however, comes with its own set of trade-offs. The FEM framework is built upon a rigorous mathematical foundation with well-established convergence guarantees. The ANN approach, while powerful, relies on a non-convex optimization process whose success is sensitive to hyperparameter tuning and initialization. Nonetheless, for this class of problem, the results demonstrate that a carefully configured physics-informed neural network can serve as a highly accurate and parameter-efficient alternative to traditional numerical methods.

\section{Conclusion}\label{conclusion}

In this investigation, two distinct numerical paradigms—the classical continuous Galerkin FEM and a modern PINNs approach—were developed and critically evaluated for the solution of the nonlinear one-dimensional Electrohydrodynamic (EHD) boundary value problem. The PINN methodology represents a significant departure from traditional discretization techniques, reformulating the differential equation as a constrained optimization problem. By embedding the governing physical laws and boundary conditions directly into a composite loss function, the neural network learns the solution without requiring a priori knowledge of the exact analytical form, a key advantage in practical engineering and scientific applications.

A primary contribution of this work is the exhaustive hyperparameter sensitivity analysis conducted within the PINN framework. The systematic investigation into the impact of data point density, activation function selection, network depth, and width on solution accuracy provides a foundational guide for applying deep learning to this class of problem. To the best of our knowledge, such a detailed comparative study for the EHD model is novel in the existing literature. The implementations of both methods were carried out using state-of-the-art software libraries—\textsf{Keras} with a \textsf{Tensorflow} backend for the PINN and the \textsf{deal.II} library for the FEM—and were rigorously validated against benchmark problems.

The comparative analysis revealed that the optimized PINN architecture achieved a marginally superior accuracy ($L_2$ error of $7.27 \times 10^{-4}$) compared to a high-resolution FEM simulation ($L_2$ error of $7.41 \times 10^{-4}$), despite utilizing significantly fewer degrees of freedom (305 trainable parameters vs. 1025 nodal values). This highlights the remarkable parameter efficiency of the deep learning approach, which constructs a global, continuous approximation of the solution. While FEM offers a mathematically rigorous framework with guaranteed convergence properties, the PINN methodology presents a compelling, mesh-free alternative that demonstrates immense potential for computational efficiency.

The promising results presented herein establish a strong case for the utility of PINNs in solving complex, nonlinear differential equations. Future work should focus on extending this methodology to higher-dimensional and more intricate physical systems. A thorough investigation into the scalability, generalizability, and theoretical convergence properties of the PINN approach for a broader class of partial differential equations will be crucial for establishing it as a robust and widely applicable tool in the scientific computing landscape.

\bibliographystyle{plain}
\bibliography{mlreferences}

\end{document}